\documentclass[15pt]{article}
\usepackage{amsfonts}

\usepackage{latexsym,wrapfig,picinpar}
\usepackage{latexsym,bm}
\usepackage{amssymb}
\usepackage{graphicx}
\usepackage{mathrsfs}
\usepackage{epsfig,latexsym}
\usepackage{amsmath}
\usepackage{amssymb}
\usepackage{fancyref}
\usepackage{epsfig}
\usepackage{dcolumn}
\usepackage{color}
\usepackage{pstricks}
\usepackage{epstopdf}
\usepackage{subfigure}
\usepackage{multicol}

\newtheorem{thm}{Theorem}[section]
\newtheorem{tm}[thm]{Theorem}
\newtheorem{lem}[thm]{Lemma}

\newtheorem{rem}[thm]{Remark}

\newtheorem{pppp}{Proof}

\newcommand{\qed}{\hspace{1em}\mbox{\raisebox{0.65ex}{\fbox{}}}}

\numberwithin{equation}{section}

\DeclareMathOperator*{\argmin}{argmin}

\allowdisplaybreaks

\begin{document}\large
\title{\bf Heterogeneous Multiscale Method for elliptic problem without scale separation}
\date{\empty}

\author{Tao Yu\footnotemark[2] ,\   Xingye Yue\footnotemark[3] \ and Changjuan Zhang\footnotemark[4]}
\maketitle

\renewcommand{\thefootnote}{\fnsymbol{footnote}}

\footnotetext[2]{
Department of Mathematics and Physics, Jinggangshan University, Ji'an 343009, P.R. China. E-mail: yutao@jgsu.edu.cn }

\footnotetext[3]{
Corresponding author. School of Mathematical Sciences and Center for Financial Engineering, Soochow University, Suzhou 215006,
P.R. China. E-mail: xyyue@suda.edu.cn. }
\footnotetext[4]{
School of Mathematical Sciences, Soochow University, Suzhou 215006,
P.R. China. E-mail: zhangchangjuan@hotmail.com. }

\date{}
\maketitle
\begin{abstract}\large
 This paper shows that the Heterogeneous Multiscale Method can be applied to elliptic problem without scale separation. The Localized Orthogonal Method is a special case of the Heterogeneous Multiscale Method.
 \end{abstract}

\section{Introduction}
\label{Intro} \setcounter{equation}{0}
Consider the classical multiscale elliptic problem
\begin{eqnarray}\label{origin-problem}
\left\{
\begin{array}{ll}\displaystyle
-\nabla\cdot(a^{\varepsilon}(x)\nabla
u(x))
=f(x),\;\;&x\in\Omega\subset\mathbb{R}^d,\\\displaystyle
u=0,&x\in\partial\Omega,\end{array}\right. 
\end{eqnarray}
where $\varepsilon\ll 1$ is a small parameter that represents the
multiscale nature in the problem. Several multiscale methods have been developed for the numerical solution of this multiscale elliptic problem. The most well known method is the Heterogeneous Multiscale Method (HMM), developed in \cite{E2003}. The method is applied into elliptic homogenization problem in \cite{E2005} and parabolic homogenization problem in \cite{Ming2006}. HMM is a general framework for designing multiscale algorithms. It consists of two components: selection of a macroscopic solver and estimating the missing macroscale data by solving locally the fine scale problem. The convergence analysis usually assumes certain periodicity and scale separation. In this paper, we deduce that HMM can be applied to multiscale elliptic problem without scale separation, and the optimal error estimate can also be obtained for this case.

Recently, a Localized Orthogonal Decomposition (LOD) method was introduced for elliptic multiscale problems in \cite{Malqvist-LOD-elliptic}, fitting into the general framework of the Variational Multiscale Method (VMS)\cite{Hughes-VMS}. The method constructs local generalized finite element basis which are exponential decayed. The analysis does not rely on high regularity of the solution or scale separation in the coefficient. The method is also applied to parabolic problems \cite{Malqvist-LOD-parabolic}. Peterseim review the VMS methods for Linear multiscale PDEs in \cite{Peterseim-VMS-decay} and give a new simply proof for the exponential decay property. In this paper, we show that the LOD method is a special case of HMM. Thus, the error estimate for HMM can be obtained from the LOD method directly.

The outline of the paper is as follows. In Section 2, we review the original definition for the Heterogeneous Multiscale Method. The localization technology is developed in Section 3.

\section{Heterogeneous Multiscale Method}
\label{} \setcounter{equation}{0}
In this section, we review the Heterogeneous Multiscale Method, first introduced in \cite{E2003}. Consider a microscopic variational problem:
\begin{equation}\label{micro-problem}
\min_{v\in V}J(v), \quad J(v)=\frac{1}{2} a(v,v)-(f,v),
\end{equation}
where $V=H_0^1(\Omega)$ and
\begin{equation}\label{micro-a}
a(u,v)=\int_{\Omega}a^{\varepsilon}\nabla u \cdot \nabla v dx,\quad (f,v)=\int_{\Omega}f v dx.
\end{equation}
The exact solution of the microscopic variational problem (\ref{micro-problem}) is
\begin{equation}\label{exact-solution}
u=\argmin_{v\in V}J(v).
\end{equation}
Let $Q:V \rightarrow V_h$ be an appropriately chosen compression operator. Here $V_h$ can be the standard $P^1$ finite element space. For each macro function $v_h\in V_h$,  define the macroscopic energy
\begin{equation}\label{macro-ennergy}
J_h(v_h)=\min_{v\in V,\;Q(v)=v_h}J(v),
\end{equation}
then the macroscopic variational problem reads
\begin{equation}\label{macro-problem}
\min_{v_h\in V_h}J_h(v_h) \triangleq \min_{v_h\in V_h}\left(\min_{v\in V,\;Q(v)=v_h}J(v)\right).
\end{equation}
Define the reconstruction operator $
R: V_h \rightarrow V$,
according to the compression operator $Q$ by
\begin{equation}\label{R}
R(w_h)=\argmin_{v\in V,\;Q(v)=w_h}J(v).
\end{equation}
The macroscopic variational problem (\ref{macro-problem}) can be rewritten as
\begin{equation}
\min_{v_h\in V}J(R(v_h)).
\end{equation}
The corresponding solution can be written as
\begin{equation}\label{uh}
u_h=\argmin_{v_h\in V}J(R(v_h)).
\end{equation}
Obviously $R(u_h)=u$.
Different choice of the compression operator $Q$ deduced different multiscale method. In this paper, we choose $Q=P_0$ be the $L_2-$projection, defined by
\begin{equation}
(P_0 v,w)=(v,w),\quad \forall w\in V_h.
\end{equation}
For this choice of $Q$, the microscopic space $V$ can be decomposed into two parts
\begin{equation}\label{space-decomposition}
V=V_h\oplus V_f,
\end{equation}
where $V_f={V_h}^{\perp}$ in the sense of $L_2-$projection $P_0$.

The following Lemma helps to understanding the reconstruction operator $
R$.
\begin{lem}\label{lagrange}
Let $u\in V$ be the solution of (\ref{R}), i.e.,
\begin{equation}
u=\underset{v\in V,\;P_0 v=v_h}{\argmin}J(v),
\end{equation}
then $P_0u=v_h$ and
\begin{equation}\label{RaVf}
a(u,w)=(f,w),\ \forall w\in {V_f}.
\end{equation}
\end{lem}
\begin{pppp}
Rewrite (\ref{RaVf}) as an minimal value problem without constrained conditions
\begin{equation}
(u,P_0\mu)=\argmin_{v\in {V},\;\mu(x)\in L_2(\Omega)}\left(J(v)+<P_0\mu,P_0v-v_h>\right),
\end{equation}
where $P_0\mu$ is the Lagrange multiplier.
Define two-variables function
\begin{equation}
G(\lambda,\epsilon)=J(u+\lambda w)+<P_0\mu+\epsilon \gamma,P_0(u+\lambda w)-v_h>,
\end{equation}
where $\lambda,\epsilon\in \mathbb{R}$, $w\in V$ and $\gamma\in L^2(\Omega)$. When $\lambda=0$ and $\epsilon=0$, function $G(\lambda,\epsilon)$ has a minimizer. Hence
\begin{equation}
\left.\frac{\partial G}{\partial \lambda}\right|_{\lambda=0,\epsilon=0}=0,
\left.\frac{\partial G}{\partial \epsilon}\right|_{\lambda=0,\epsilon=0}=0.
\end{equation}
That is to say
\begin{equation}\label{glam}
\left.\frac{\partial G}{\partial \lambda}\right|_{\lambda=0,\epsilon=0}=a(u,w)-(f,w)+(P_0\mu,P_0w)=0,\quad \forall w\in V,
\end{equation}
\begin{equation}\label{geps}
\left.\frac{\partial G}{\partial \epsilon}\right|_{\lambda=0,\epsilon=0}=(\gamma,P_0u-v_h)=0,\quad \forall \gamma\in L^2(\Omega).
\end{equation}
From \eqref{geps}, we get
\begin{equation}
P_0u=v_h.
\end{equation}
Let $w=w_f\in{V_f}\subset V$, then $P_0 w=0$ by the definition of $V_f$. Thus $(P_0u,P_0w)=0$ and
\begin{equation}\label{u-equation}
a(u,w)=(f,w),\ \forall w\in {V_f},
\end{equation}
from (\ref{glam}).
Denote $u=P_0 u+\Phi=v_h+\Phi$, where $\Phi\in{V_f}$. Brought it into (\ref {u-equation}),
\begin{equation}
a(\Phi,w)=(f,w)-a(v_h,w),\ \forall w\in {V_f}.
\end{equation}
Through this, we can detest $\Phi$ and then  determine the solution $u$. After the solution $u$ has determined and from \eqref{glam}, \begin{equation}(P_0\mu,P_0w)=<<F,w>>\triangleq-a(u,w)+(f,w),
\end{equation}
which can determine the Lagrange multiplier $P_0\mu$.
\end{pppp}

Define an approximation of the reconstruction operator $R$ by
\begin{equation}\label{Rh}
R_h(v_h)=\argmin_{v\in {V},P_0(v)=v_h}J_0(v),
\end{equation}
where $J_0(v)=\frac{1}{2}a(v,v)$. Then the corresponding approximation solution $u_h^{ms}$ is defined by
\begin{equation}\label{uhms}
u_h^{ms}=R_h(\tilde{u}_h)=\argmin_{v_h\in R_h({V_h})}J(v_h),
\end{equation}
where
\begin{equation}\label{uth}
\tilde{u}_h=\argmin_{v_h\in {V_h}}J(R_h(v_h)).
\end{equation}

The following Lemma shows the error between $R(\tilde{u}_h)$ and $R_h(\tilde{u}_h)$.
\begin{lem}\label{rf}
$R$ and $R_h$ are defined by (\ref{R}) and (\ref{Rh}), respectively. Then
 $$R(\tilde{u}_h)=R_h(\tilde{u}_h)+R_f(f),$$
 where $R_f(f)=\underset{v\in {V_f}}{\argmin}J(v)$.
\end{lem}
\begin{pppp}
From Lemma \ref{lagrange}, $R(\tilde{u}_h)$ satisfies
\begin{equation}\label{aR}
a(R(\tilde{u}_h),v)=(f,v),\ \forall v\in {V_f}.
\end{equation}
Similarly, $R_h(\tilde{u}_h)$ satisfies
\begin{equation}\label{aRh}
a(R_h(\tilde{u}_h),v)=0,\ \forall v\in {V_f}.
\end{equation}
Denote $R_f(f)=R(\tilde{u}_h)-R_h(\tilde{u}_h)$, then
\begin{equation}
a(R_f(f),v)=a(R(\tilde{u}_h)-R_h(\tilde{u}_h),v)=(f,v),\ \forall v\in\ {V_f},
\end{equation}
by subtracting (\ref{aR}) and (\ref{aRh}).
That is to say $R_f(f)=\underset{v\in {V_f}}{\argmin}J(v)$.
\end{pppp}

\begin{lem}\label{equalu}
$u_h$ and $\tilde{u}_h$ are defined by (\ref{uh}) and (\ref{uth}), respectively. Then they are the same solution, i.e.
\begin{equation}
\tilde{u}_h=u_h.
\end{equation}
\end{lem}
\begin{pppp}
Denote $R_f(v_h)\triangleq R({v_h})-R_h({v_h}),\;\forall v_h\in V_h$, then
$$P_0(R_f(v_h))=P_0(R({v_h}))-P_0(R_h({v_h}))=v_h-v_h=0$$
by the definition of $R$ and $R_h$, i.e. $R_f(v_h)\in {V_f}$.
From Lemma \ref{rf}, the remainder term $R_f(v_h)=\underset{v\in {V_f}}{\argmin}J(v)$ does not relied on $v_h$. Abbreviate $R_f(v_h)$ by $R_f$ and from the definition of energy $J$, then
\begin{eqnarray}
J(R(v_h))\!\!\!\!&=&\!\!\!\!\frac{1}{2}a(R(v_h),R(v_h))-(f,R(v_h))\nonumber\\
\!\!\!\!&=&\!\!\!\!\frac{1}{2}a(R_h(v_h)+R_f,R_h(v_h)+R_f)-(f,R_h(v_h)+R_f) \nonumber\\\!\!\!\!&=&\!\!\!\!\frac{1}{2}a(R_h(v_h),R_h(v_h))- (f,R_h(v_h))+\frac{1}{2}a(R_f,R_f)-(f,R_f)+a(R_h(v_h),R_f)\nonumber\\ \!\!\!\!&=&\!\!\!\!J(R_h(v_h))+J(R_f),\nonumber
\end{eqnarray}
by $a(R_h(v_h),w)=0,\;\forall w\in V_f$ and $R_f\in V_f$. As $R_f$ does not relied on $v_h$, $J(R_f)$ also does not relied on $v_h$. Thus
$$\argmin_{v_h\in {V_h}}J(R(v_h))=\argmin_{v_h\in {V_h}}J(R_h(v_h)),$$
that is to say $u_h=\tilde{u}_h$.
\end{pppp}
Define the energy norm $|||\cdot|||:=|||\cdot|||_{\Omega}:=||a^{\varepsilon}\nabla\cdot||_{L^2(\Omega)}$ on space $V$. The error estimate between the exact solution $u$ and the approximation solution $u_h^{ms}$ is deduced by the following Theorem.
\begin{tm}\label{error-estimate}
Let $u\in V$ be the exact solution of (\ref{micro-problem}) and $u_h^{ms}$ be the approximation solution defined in (\ref{uhms}). Then it holds that
$$|||u-u_h^{ms}|||\leq Ch||f||_0$$
with positive constant $C$ that does not depend on $h$ and $\varepsilon$.
\end{tm}
\begin{pppp}
Due to Lemma \ref{rf} and Lemma \ref{equalu}, it holds that
$$u-u_h^{ms}=R(u_h)-R_h(\tilde{u_h})=R(\tilde{u_h})-R_h(\tilde{u_h})=R_f(f)$$ by the definition (\ref{micro-problem}) and (\ref{uhms}). The remainder term $R_f(f)$ satisfies $a(R_f(f),v)=(f,v),\ \forall v\in {V_f}$ from Lemma \ref{rf}. The application of $L^2$-projection error estimate and H\"{o}lder inequality yield
\begin{equation}
\begin{split}
|||u-u_h^{ms}|||^2&=|||R_f(f)|||^2=a(R_f(f),R_f(f))\\
&=(f,R_f(f))\leq ||f||_0||R_f(f)||_0\\
&=||f||_0||R_f(f)-P_0(R_f(f))||_0\\
&\leq Ch||f||_0|||R_f(f)|||,
\end{split}
\end{equation}
since $P_0(R_f(f))=0$. This concludes the proof.
\end{pppp}
\begin{rem}
For every $v\in {V_f}$, it holds
$$||v||_0=||v-P_0v||_0\leq Ch||v||_1,$$
where $C$ be a positive constant.
\end{rem}

\section{Localization}
\label{} \setcounter{equation}{0}
From Theorem \ref{error-estimate}, we know that $u_h^{ms}$ is a good approximation of the exact solution $u$. But $u_h^{ms}$ is solved in the whole domain $\Omega$, which is too expensive to compute. In this section, we deduce that we can compute it locally, and then the same error estimate can be obtained.

For every vertex $z\in \mathcal{N}$, the vertex set of $V_h$, denote the corresponding nodal basis function $\lambda_z$, satisfying
\begin{equation}
\lambda_z(z)=1\ and\ \lambda_z(x)=1\ for\ all\ x \neq z\in \mathcal{N}.
\end{equation}
Define corresponding basis function $R_h(\lambda_z)=\underset{v\in{V},P_0(v)=\lambda_z}{\argmin}J_0(v)$, which satisfies
\begin{equation}\label{Rhlz}
a(R_h(\lambda_z),w)=0,\ \forall w\in V_f,
\end{equation}
as Lemma \ref{lagrange} similarly.
From the space decomposition (\ref{space-decomposition}), $R_h(\lambda_z)$ can be decomposed into
\begin{equation}\label{Rhlzdc}
R_h(\lambda_z)=P_0(R_h(\lambda_z))+\phi_z=\lambda_z+\phi_z,
\end{equation}
where the corrector $\phi_z\in V_f$ satisfies
\begin{equation}\label{phiz}
a(\phi_z,w)=-a(\lambda_z,w),\ \forall w\in V_f,
\end{equation}
from (\ref{Rhlz}).
These corresponding basis functions from basises of $R_h(V_h)$, that is to say,
\begin{equation}
R_h({V_h})=span\{R_h(\lambda_z)|z\in \mathcal{N}\}=span\{\lambda_z+\phi_z|z\in \mathcal{N}\}.
\end{equation}
Thus, the multiscale method (\ref{uhms}) can be rewritten as: find $u_h^{ms}\in R_h({V_h})$, such that
\begin{equation}
a(u_h^{ms},v)=(f,v),\ \forall v\in R_h({V_h}).
\end{equation}
The following Lemma shows that the corrector $\phi_z$ is exponential decayed, which is proved in \cite{Peterseim-VMS-decay}.
\begin{lem}\label{}
There exists a constant $c>0$ independent of $h$ and $R$ such that
$$|||\phi_z|||_{\Omega\setminus B_R(z)}\leq \exp\left(-c\frac{R}{h}\right)|||\phi_z|||_{\Omega},$$
where $B_R(z)$ denotes the ball of radius $R>0$ centered at $z$.
\end{lem}

Denote nodal patch $\omega_{z,k},\ z\in \mathcal{N}$ by
\begin{equation}
\begin{split}
&\omega_{z,1}=supp(\lambda_z)=int\left(\cup\{K\in \mathcal{T}_h|z\in K\}\right),\\
&\omega_{z,k}=int\left(\cup\{K\in \mathcal{T}_h|K\cap\bar{\omega}_{z,k-1}\neq\varnothing\}\right),\ k=2,3,4,\cdots.
\end{split}
\end{equation}
Define localized basis function
\begin{equation}
\tilde{R}_h(\lambda_z)=\argmin_{v\in H_0^1(\omega_{z,k}),P_0(v)=\lambda_z}J_0(v),
\end{equation}
then the corresponding localized space $\tilde{R}_h({V_h})$ is obtained by
\begin{equation}
\tilde{R}_h({V_h})=span\{\tilde{R}_h(\lambda_z)|z\in \mathcal{N}\}.
\end{equation}
Then the localized approximation solution is defined by
\begin{equation}\label{uhkms}
u_{h,k}^{ms}=\argmin_{v_h\in \tilde{R}_h({V_h})}J(v_h).
\end{equation}
That is to say, the localized multiscale method reads: find $u_{h,k}^{ms}\in \tilde{R}_h({V_h}) $, such that
\begin{equation}\label{}
a(u_{h,k}^{ms},v)=(f,v),\ \forall v\in \tilde{R}_h({V_h}).
\end{equation}
Next, we shall show the expression of $\tilde{R}_h(\lambda_z)$. Similarly as (\ref{Rhlz}) and (\ref{Rhlzdc}), it's easily shown that
\begin{equation}\label{tRhlzdc}
\tilde{R}_h(\lambda_z)=\lambda_z+\phi_{z,k},
\end{equation}
where the localized corrector $\phi_{z,k}\in V_f\cap H_0^1(\omega_{z,k})$ satisfies
\begin{equation}\label{phiz}
a(\phi_{z,k},w)=-a(\lambda_z,w),\ \forall w\in V_f\cap H_0^1(\omega_{z,k}).
\end{equation}

The error estimate between $u$ and $u_{h,k}^{ms}$ is shown in the following Theorem, which can be found in \cite{Peterseim-VMS-decay}.
\begin{tm}\label{error-estimate-local}
Let $u\in V$ be the exact solution of (\ref{micro-problem}) and $u_{h,k}^{ms}$ be the localized approximation solution defined in (\ref{uhkms}). Then it holds that
$$|||u-u_{h,k}^{ms}|||\leq Ch||f||_0$$
with moderate choice $k$ and positive constant $C$ that does not depend on $h$ and $\varepsilon$.
\end{tm}

\section{Acknowledgement}
This work was supported by the Natural Science Foundation of China (No. 62241203) and Science and Technology Research Project of Jiangxi Provincial Department of Education(No. GJJ211027).


\end{document}